\documentclass[12pt]{amsproc}
\usepackage{amsmath,amssymb,amscd,stmaryrd}
\usepackage[mathscr]{euscript}

\newcommand{\sB}{\mathscr B}
\newcommand{\sC}{\mathscr C}

\newcommand{\sE}{\mathscr E}
\newcommand{\sF}{\mathscr F}
\newcommand{\sH}{\mathscr H}
\newcommand{\sJ}{\mathscr J}
\newcommand{\sK}{\mathscr K}
\newcommand{\sM}{\mathscr M}
\newcommand{\sP}{\mathscr P}
\newcommand{\sR}{\mathscr R}

\newcommand{\sU}{\mathscr U}
\newcommand{\sX}{\mathscr X}

\newcommand{\bC}{\mathbb C}

\renewcommand{\d}{\delta}

\begin{document}

\title[Some Results and a $K$-theory Problem]{Some Results and a $K$-Theory Problem About Threshold Commutants mod Normed Ideals}
\author{Dan-Virgil Voiculescu}
\address{D.V. Voiculescu \\ Department of Mathematics \\ University of California at Berkeley \\ Berkeley, CA\ \ 94720-3840}
\thanks{Research supported in part by NSF Grant DMS-1665534.}
\dedicatory{Dedicated to Ciprian Foias on the occasion of his 85th birthday.}

\begin{abstract}
We extend to the case of a threshold ideal our result with J.~Bourgain about the essential centre of the commutant mod a diagonalization ideal for a $n$-tuple of commuting Hermitian operators. We also compute the $K_0$-group of the commutant mod trace-class of a unitary operator with spectrum equal to its essential spectrum. We present the problem of computing the $K_1$-group for a commutant mod trace-class in its simplest case.
\end{abstract}
\maketitle

\section{Introduction}
\label{sec1}

The commutant modulo a normed ideal of a $n$-tuple of Hermitian operators is a Banach algebra with involution, which is not a $C^*$-algebra. However, if the normed ideal is sufficiently large, then the quotient of this Banach algebra by its ideal of compact operators is a $C^*$-algebra. An interesting situation arises for certain threshold normed ideals, when the quotient by the compact ideal is not yet a $C^*$-algebra, but is isomorphic as a Banach algebra with involution to a $C^*$-algebra, that is the quotient norm is equivalent to a $C^*$-norm. In technical terms of the quasicentral modulus $k_{\sJ}(\tau)$ of the $n$-tuple $\tau$ with respect to the normed ideal $\sJ$ the first situation occurs when $k_{\sJ}(\tau)$ vanishes, while the situation of the threshold ideal is when $k_{\sJ}(\tau)$ is finite but non-zero (for details and references see our survey paper \cite{12}). When $n=1$, the case of one Hermitian operator $T$, the threshold situation occurs when the normed ideal is the trace-class $\sC_1$ and $T$ has absolutely continuous spectrum of bounded multiplicity (a sufficient condition). Note that in this setting the Kato--Rosenblum theorem and trace-class scattering theory can be used.

The present paper continues the study of threshold commutants mod normed ideals for $n$-tuples of commuting Hermitian operators. We deal with three questions. One result is about the centre of the quotient by the compact ideal, which we show, roughly, is generated by the classes of components of $\tau$. This extends to the threshold, the result of \cite{3} in the case of vanishing $k_{\sJ}(\tau)$. We should say that this doesn't mean that the quotient algebra is isomorphic to the $C^*$-algebra of a continuous field of $C^*$-algebras over the essential spectrum of $\tau$ (this is pointed out in the discussion in section~\ref{sec5}). A second result is about the $K_0$-group of the commutant mod trace-class of a unitary operator. This is an extension of the result for one Hermitian operator in \cite{11}, with a new feature that roughly the $K_0$-class views the Hardy-subspace as absolutely continuous with multiplicity $1/2$. The third question we bring up is the computation of the $K_1$-group in the simplest case. The remarks we make about this problem suggest that techniques from trace-class scattering theory may be relevant.

The paper, besides the introduction and references has four more sections. Section~\ref{sec2} contains preliminaries. Section~\ref{sec3} proves the result about the essential centre. Section~\ref{sec4} deals with $K_0$ in the case of a unitary operator. Section~\ref{sec5} is a discussion around $K_1$ in the simplest case.

\section{Preliminaries}
\label{sec2}

Throughout this paper $\sH$ will be a complex separable infinite dimensional Hilbert space and $\sB(\sH)$, $\sK(\sH)$, $\sB/\sK(\sH)$, or simply $\sB$, $\sK$, $\sB/\sK$ will denote the bounded operators, the compact operators and the Calkin algebra. The canonical homomorphism $\sB \to \sB/\sK$ will be denoted by $p$.

If $\tau = (T_j)_{1 \le j \le n}$ is a $n$-tuple of bounded Hermitian operators and $(\sJ,|\quad|_{\sJ})$ is a normed ideal (see \cite{4}, \cite{9}) we denote by
\[
\sE(\tau;\sJ) = \{X \in \sB \mid [T_j,X] \in \sJ,\ 1 \le j \le n\}
\]
the commutant mod~$\sJ$ of $\tau$, which is a Banach algebra with involution when endowed with the norm
\[
\|| X\|| = \|X\| + \max_{1 \le j \le n} |[X,T_j]|_{\sJ}.
\]
We also denote by
\[
\begin{aligned}
\sK(\tau;\sJ) &= \sE(\tau;\sJ) \cap \sK \mbox{ and by} \\
\sE/\sK(\tau;\sJ) &= \sE(\tau;\sJ)/\sK(\tau;\sJ)
\end{aligned}
\]
its compact ideal and the quotient Banach algebra.

The notation also makes sense when we drop the condition that the $T_j$ be Hermitian. The algebra is still involutive when we require only that $\{T_1,\dots,T_n\} = \{T^*_1,\dots, T^*_n\}$ and also in case $n=1$ when $T_1  = U$ is a unitary operator.

Note also that $\sE/\sK(\tau;\sJ)$ is algebraically isomorphic to $p(\sE(\tau;\sJ)) \subset \sB/\sK$ and that the map
\[
\sE/\sK(\tau;\sJ) \to p(\sE(\tau;\sJ))
\]
is contractive.

The quasicentral modulus of $\tau$ with respect to $\sJ$ is the number
\[
k_{\sJ}(\tau) = \liminf_{A \in \sR^+_1} \max_{1 \le j \le n} |[A,T_j]|_{\sJ}
\]
where $\sR^+_1$ is the set of finite rank positive contractions on $\sH$ endowed with the natural order (see \cite{12}).

If $k_{\sJ}(\tau) = 0$ and the finite rank operators are dense in $\sJ$, then $\sE/\sK(\tau;\sJ)$ is a $C^*$-algebra, while if only $k_{\sJ}(\tau) < \infty$, $\sE/\sK(\tau;\sJ)$ is a Banach algebra with involution isomorphic to a $C^*$-algebra and $p(\sE(\tau;\sJ))$ is a $C^*$-subalgebra of $\sB/\sK$ (see \cite{12} and references therein).

By $\sM_n$ we denote the $n \times n$ matrices. We have the identification
\[
\sM_n\sE(\tau;\sJ) \sim \sE(\tau \otimes I_n;\sJ).
\]
Here $\tau \otimes I_n$ is viewed as acting on $\sH \otimes \bC^n \simeq \sH^n$. The above identification is up to equivalent norms, but this will not matter in the questions we will consider.

By $(\sC_p,|\,|_p)$ we denote the Schatten--von~Neumann $p$-classes.

The idempotents in $\sE(\tau;\sJ)$ will be denoted by $\sP\sE(\tau;\sJ)$ and the Hermitian idempotents will be denoted by $\sP_h\sE(\tau;\sJ)$.

If $\tau$ is a $n$-tuple of commuting normed operators we denote by $\sigma(\tau)$ the joint spectrum and by $\sigma_e(\tau) = \sigma(p(\tau))$ the joint essential spectrum.

\section{Commutation mod $\sK$}
\label{sec3}

The commutation mod~$\sK$ result we prove in this section is a stronger result which has as a corollary the essential centre result for threshold commutants mod normed ideals which extends to the threshold the result of \cite{3}.

\bigskip
\noindent
{\bf Theorem 3.1.} {\em 
Let $\tau = (T_j)_{1 \le j \le n}$ be a $n$-tuple of commuting Hermitian operators, the joint spectrum of which $K = \sigma(\tau) \subset {\mathbb R}^n$ is a perfect set. If $X$ is so that
\[
[X,\sE(\tau;\sC_1)] \subset \sK
\]
then we have $X \in C^*(\tau) + \sK$.
}

\bigskip
\noindent
{\em Proof.} 
Clearly it suffices to prove the assertion when $X = X^*$.

Remark also that $\tau$ can be replaced by another $n$-tuple of commuting Hermitian operators which has the same spectrum as $\tau$ and which is a trace-class perturbation of $\tau$. Indeed this will not change $\sE(\tau;\sC_1)$ and $C^*(\tau) + \sK$. Using for instance the adaptation of the Voiculescu theorem to normed ideals \cite{10}, we can find such a modified $\tau$ which is unitarily equivalent to $\tau \oplus \d$ where $\d$ is a given $n$-tuple of commuting Hermitian operators with $\sigma(\d) \subset K$ and which is diagonalizable in some orthonormal basis. Based on this we can assume that for every $k \in K$ the joint eigenspace $E(\tau;\{k\})\sH$ is either zero or infinite dimensional. Indeed, for this we choose $\d$ to be an infinite ampliation of the restriction of $\tau$ to the direct sum of the $E(\tau;\{k\})\sH$ where $k \in K$.

Since $(\tau)' \subset \sE(\tau;\sC_1)$ we have $[X,(\tau)'] \subset \sK$ and hence by \cite{5} it follows that $X \in (\tau)'' + \sK$. This means there is a bounded Borel function $f: K \to {\mathbb R}$ so that $X \in f(\tau) + \sK$. Thus the proof reduces to the case when $X = f(\tau)$.

Identifying ${\mathbb R}^n$ with $\ell^{\infty}(\{1,\dots,n\})$ we shall use on ${\mathbb R}^n$ the metric arising from the $\ell^{\infty}$-norm. If $\omega \subset K$ is a Borel set, let $\Delta(\omega) = \mbox{diam}(\sigma(X \mid E(\tau;\omega)\sH)$ that is the diameter of the spectrum of the restriction of $X$ to the spectral subspace of $\tau$ for $\omega$.

We shall prove that if $\mbox{diam } \omega_j \to 0$ for a sequence of Borel sets $\omega_j \subset K$ then we have $\Delta(\omega_j) \to 0$.

Suppose the contrary that $\mbox{diam } \omega_j \to 0$ and $\Delta(\omega_j) \nrightarrow 0$ for some sequence. Passing to a subsequence we can assume there is $k_0 \in K$ so that $r_j = d(\omega_j,k_0) \downarrow 0$ and that all $\Delta(\omega_j) > \varepsilon$ for some $\varepsilon > 0$.

On the other hand, observe that if $\Omega_p \uparrow \Omega$ then $\Delta(\Omega_p) \uparrow \Delta(\Omega)$. In particular, we have
\[
\begin{aligned}
&\Delta(\omega_j\backslash(B(k_0,r_n)\backslash\{k_0\})) \uparrow \Delta(\omega_j) \mbox{ and} \\
&\Delta(\omega_j\backslash B(k_0,r_n)) \uparrow \Delta(\omega_j\backslash\{k_0\} \mbox{ as $n \to \infty$,}
\end{aligned}
\]
where $B(k_0,r_n)$ denotes the ball in the $\ell^{\infty}$ metric. In particular, passing to a subsequence and replacing then $\omega_j$ by $\omega_j\backslash (B(k_0;r_{j+1})\backslash \{k_0\})$ we may assume the $\omega_j\backslash \{k_0\}$ are pairwise disjoint. Remark also that the spectrum of $X \mid E(\tau;\omega_j)\sH$ 
equals $\sigma(X \mid E(\tau;\omega_j \backslash \{k_0\})\sH$ or $\sigma(X \mid E(\tau;\omega_j\backslash \{k_0\})\sH) \cup \{f(k_0)\}$ depending on whether $E(\tau;\{k_0\})$ is $= 0$ or $\ne 0$. We consider now two cases.

\medskip
a) If $\liminf_{j \to \infty} \Delta(\omega_j\backslash\{k_0\}) = 0$, passing to a subsequence so that $\Delta(\omega_j\backslash\{k_0\}) \to 0$, we must have $k_0 \in \omega_j$ and $E(\tau;\{k_0\}) \ne 0$. Moreover since $\Delta(\omega_j) > \varepsilon > 0$, there is $\lambda_j \in \sigma(X \mid E(\tau;\omega_j\backslash \{k_0\})\sH)$ so that $|\lambda_j - f(k_0)| > \varepsilon/2$ for each $j$. Recall also that by our additional assumptions on $\tau$ we have that $E[\tau;\{k_0\}]\sH$ is infinite dimensional.

We can then choose for each $j$ a number $\mu_j \in \{\lambda_j,f(k_0)\}$ so that $|\mu_j-\mu_{j+1}| > \varepsilon/4$ for all $j \in {\mathbb N}$. We also choose an orthonormal sequence $(h_j)_{j \in {\mathbb N}}$ of vectors
\[
h_j \in E(\tau;\omega_j)\sH \cup E(\tau;\{k_0\})\sH
\]
so that
\[
\|Xh_j - \mu_jh_j\| < 10^{-j}.
\]
Here $h_j \in E(\tau;\omega_j\backslash\{k_0\})\sH$ if $\mu_j = \lambda_j$ and $h_j \in E(\tau;\{k_0\})\sH$ if $\mu_j = f(k_0)$. Finding such an orthonormal sequence is possible because the $\omega_j\backslash\{k_0\}$ are disjoint and $E(\tau;\{k_0\})\sH$ is infinite dimensional.

Consider then the shift operator on the sequence $(h_j)_{j \in {\mathbb N}}$:
\[
Y = \sum_j \langle \cdot,h_j\rangle h_{j+1}.
\]
We have that
\[
\begin{aligned}
&\|[T_k,\langle \cdot,h_j\rangle h_{j+1}]\| = \\
&= \|\langle \cdot,T_kh_j\rangle h_{j+1} - \langle \cdot,h_j\rangle T_kh_{j+1}\| \le \\
&\le 10(r_j + r_{j+1}).
\end{aligned}
\]

Passing again if necessary to a subsequence we may assume the sequence of $r_j$'s is summable and hence
\[
[Y,T_k] \in \sC_1,\ 1 \le k \le n
\]
so that $Y \in \sE(\tau;\sC_1)$.

On the other hand
\[
[X,Y] - \sum_j(\langle\cdot,h_j\rangle\mu_{j+1}h_{j+1} - \langle\cdot,\mu_jh_j\rangle h_{j+1})
\]
is compact since
\[
\begin{aligned}
&\|[X,\langle\cdot h_j\rangle h_{j+1}] - (\mu_{j+1}-\mu_j)\langle\cdot h_j,h_{j+1}\rangle\| \\
&\le \|X h_j-\mu_jk_j\| + \|Xh_{j+1} - \mu_{j+1}h_{j+1}\| \\
&\le 2 \cdot 10^{-j}.
\end{aligned}
\]
Since
\[
\sum_j(\mu_{j+1}-\mu_j)\langle \cdot h_j,h_{j+1}\rangle
\]
is not compact, because $|\mu_{j+1}-\mu_j| \ge \varepsilon/4$, we have arrived at a contradiction having $[Y,X] \notin \sK$, while $Y \in \sE(\tau;\sC_1)$.

\medskip
b) The other case being when $\liminf_{j \to \infty} \Delta(\omega_j\backslash\{k_0\}) > 0$ which means $\inf_{j \in {\mathbb N}} \Delta(\omega_j\backslash\{k_0\}) > 0$, we can find $\alpha_j,\beta_j \in \sigma(X \mid E(\tau;\omega_j\backslash\{k_0\})\sH)$ so that $|\alpha_j-\beta_j| > \eta$ for some $\eta > 0$ for all $j \in {\mathbb N}$. Then we can choose $\mu_j \in \{\alpha_j,\beta_j\}$ so that $|\mu_j-\mu_{j-1}| > \eta/2$. We can also find unit vectors $h_j \in E(\tau;\omega_j\backslash\{k_0\})\sH$ so that $\|Xh_j-\mu_jh_j\| < 10^{-j}$. The $\omega_j\backslash\{k_0\}$ being disjoint, this sequence will be orthonormal and we can proceed along the same lines as in case a). We consider the shift operator
\[
Y = \sum_j \langle \cdot,h_j\rangle h_{j+1}
\]
and like in case a) we arrive at the contradiction that $[Y,T_k] \in \sC_1$, $1 \le k \le n$, while $[Y,X] \notin \sK$.

Concluding we have proved that $\mbox{diam}(\omega_j) \to 0$ implies $\Delta(\omega_j) \to 0$. We shall use this to reach the desired conclusion that $X = f(\tau) \in C^*(\tau)$.

Thus for every $\varepsilon > 0$ there is $\varphi(\varepsilon) > 0$ so that
\[
\mbox{diam}(\omega) < \varphi(\varepsilon) \Rightarrow \Delta(\omega) < \varepsilon.
\]
Then given an integer $m > 0$ we construct a continuous function $g_m: K \to {\mathbb R}$ as follows.

Let $100 \delta = \varphi(2^{-m})$ and let $(\omega_{m,p})_{1 \le p \le N}$ be the partition of $K$ into the Borel sets which are the nonempty $K \cap (x + \delta[0,1)^n)$ where $x \in \delta {\mathbb Z}^n$. Observe that $\overline{\omega}_{m,p} \cap \overline{\omega}_{m,q} \ne \emptyset \Rightarrow d(\omega_{m,p},\omega_{m,q}) \le 2\delta$. Let further $(\theta_p)_{1 \le p \le N}$ be a continuous partition of unity $\theta_p: K \to [0,1]$ so that $d(\omega_{m,p},\mbox{supp } \theta_p) \le 3\delta$. For each $1 \le p \le N$ let also $s_p \in \sigma(X \mid E(\tau;\omega_{m,p})\sH)$. We define $g_m = \sum_{1 \le p \le N} s_p\theta_p$. Then, if $\theta_q\mid \omega_{m,p} \ne 0$ we must have $d(\omega_{m,p},\omega_{m,q}) < 10\delta$ and hence $|s_p-s_q| \le 2^{-m}$. This gives for $k \in \omega_{m,p}$
\[
|g_m(k) - s_p| \le \sum_q \theta_q(k)|s_p-s_q| \le 2^{-m}.
\]
Since
\[
\left\| \sum_{1 \le p \le N} s_pE(\tau;\omega_{m,p}) - X\right\| \le 2^{-m},
\]
we infer that $\|g_m(\tau) - X\| \le 2^{-m+1}$. On the other hand $g_m(\tau) \in C^*(\tau)$ so that $X \in C^*(\tau)$.\hfill\qed

\bigskip
\noindent
{\bf Theorem 3.2.} {\em 
Let $\tau$ be a $n$-tuple of commuting Hermitian operators the spectrum $\sigma(\tau)$ of which is a perfect set and let $\sJ$ be a normed ideal in which the finite rank operators are dense so that $k_{\sJ}(\tau) < \infty$. Then the centre of the $C^*$-algebra $p(\sE(\tau;\sJ))$ is generated by $p(\tau)$, so that its spectrum is $\sigma(\tau)$.
}

\bigskip
\noindent
{\em Proof.} 
The fact that $p(\sE(\tau;\sJ))$ is a $C^*$-algebra because of $k_{\sJ}(\tau) < \infty$ and the fact that $[T_k,\sE(\tau;\sJ)] \subset \sK$, $1 \le k \le n$, implies that the $C^*$-algebra generated by $p(\tau)$ is contained in the centre of $p(\sE(\tau;\sJ))$. The opposite inclusion follows from $X \in \sE(\tau;\sJ))$. The opposite inclusion follows from $X \in \sE(\tau;\sJ) + \sK$ and $[X,\sE(\tau;\sJ)] \subset \sK$ then $[X,\sE(\tau;\sC_1)] \subset \sK$ and Theorem~$3.1$ gives $X \in C^*(\tau) + \sK$. Note also that since $\sigma(\tau)$ is a perfect set we have $\sigma(\tau) = \sigma(p(\tau))$.\hfill\qed

\section{The case of a unitary operator}
\label{sec4}

In this section we compute the ordered group $K_0(\sE(U;\sC_1))$ where $U$ is a unitary operator with $\sigma(U) = \sigma_e(U)$. Our previous result in \cite{11} for a Hermitian operator is a particular case. Indeed if $T = T^*$ and $U = \exp(i \alpha T)$ where $\alpha > 0$, $\alpha\|T\| < \pi$ then $\sE(U;\sC_1) = \sE(T;\sC_1)$. Conversely if $\sigma(U) \ne {\mathbb T}$ then there is $T = T^*$ so that $\sE(U;\sC_1) = \sE(U,U^*;\sC_1)$. Thus the novelty will appear in case $\sigma(U) = {\mathbb T}$.

Remark that $\sE(U;\sC_!) = \sE(U,U^*;\sC_1)$ and also the norms are equal. The equivalence relation $P \sim Q$ in $\sP\sE(U;\sC_1)$, $P = XY$, $Q = YX$ with $X,Y \in \sE(U;\sC_1)$, like in \cite{11} is easily seen to be equivalent to
\[
P_{P\sH} = V^*V,\ P_{Q\sH} = VV^*
\]
where $P_{\sX}$ is the orthogonal projection onto $\sX$ and $V$ is a partial isometry in $\sE(U;\sC_1)$. Moreover $P \in \sP\sE(U;\sC_1)$ implies $P_{P\sH} \in \sE(U;\sC_1)$. We also have $P \sim P_{P\sH}$. Thus we can work with Hermitian projections.

Also, as pointed out in \cite{11}, if $P,Q \in \sP_h(\sE(U;\sC_1))$ then $P \sim Q$ iff there is a partial isometry $V \in \sB(\sH)$ so that $VV^* = P$, $V^*V = Q$ and $VQUQ - PUPV \in \sC_1$ and this then also implies that $V \in \sE(U;\sC_1)$.

The absolutely continuous and singular subspaces for $U$ will be denoted by $\sH_{ac}(U)$ and $\sH_{\sin g}(U)$ and the corresponding projections $E_{ac}(U)$ and $E_{\sin g}(U)$. The multiplicity function of $U|\sH_{ac}(U)$ will be denoted by $mac(U)$ and is a Borel function defined on ${\mathbb T}$ up to Haar null sets and taking values in $\{0,1,2,\dots,\infty\}$. We shall also denote by $\omega(U)$ the set $mac^{-1}(\infty)$, defined up to a null-set. When dealing with matrices we pass to $U \otimes I_n$ and $\sH^n$ as outlined in the preliminaries.

It will be convenient to consider an extension of the multiplicity function to deal with unitaries mod trace-class, that is the set
\[
U\sC_1(\sH) = \{F \in \sB(\sH) \mid F^*F - I \in \sC_1,FF^* - I \in \sC_1\}.
\]
Remark that if $F \in \sU\sC_1(\sH)$ then $F$ is a Fredholm operator and if $\mbox{ind } F = 0$ then there is a unitary operator $F_1$ so that $F - F_1 \in \sC_1$ and the multiplicity function of the absolutely continuous part of $F_1$ does not depend on the choice of $F_1$. Indeed if $F_2$ is another unitary operator so that $F-F_1 \in \sC_1$ then $F_1 - F_2 \in \sC_1$ and by \cite{1} we have $F_1 \mid \sH_{ac}(F_1)$ and $F_2 \mid \sH_{ac}(F_2)$ are unitarily equivalent. The existence of $F_1$ can be seen as follows, first replace $F$ by ${\tilde F}$ so that ${\tilde F}$ is invertible and $F - {\tilde F}$ is finite rank and then $F_1 = {\tilde F}({\tilde F}^*{\tilde F})^{-1/2}$ will do. Thus on the set
\[
\Omega = \{F \in \sU\sC_1(\sH) \mid \mbox{ind } F = 0\}
\]
there is a well-defined essential absolutely continuous multiplicity function $emac(F)$ so that if
$F_1,F_2 \in \Omega$ and $F_1 - F_2 \in \sC_1$ then $emac(F_1) = emac(F_1)$ Haar a.e.\ and if $F \in \Omega$ is unitary, then $mac(F) = emac(F)$ Haar a.e.  Remark also that $emac(F_1 \oplus F_2) = emac(F_1) + emac(F_2)$. To extend $emac$ to all of $\sU\sC_1(\sH)$ let $\sJ$ be a conjugate-linear antiunitary involution of $\sH$ and define
\[
emac(F) = \frac {1}{2} emac(F \oplus \sJ F^*\sJ)
\]
for any $F \in \sU\sC_1(\sH)$. Indeed $\mbox{ind}(F \oplus \sJ F^*\sJ) = 0$ so the right-hand side has been defined. Clearly the definition does not depend on the choice of $\sJ$ and on the choice of $F \mbox{ mod } \sC_1$. Remark also that if $F \in \sU\sC_1(\sH)$ and $pol(F)$ is the partial isometry from $F^*\sH$ to $F\sH$ then $F - pol(F) \in \sC_1$ so that $emac(F) = emac(pol(F))$. We can now proceed as follows to see what $emac(F)$ is. We can extend $pol(F)$ to $V$ which is an isometry or co-isometry depending on whether the index is $\le 0$ or $\ge 0$. Then $emac(pol(F)) = emac(V)$. The Wold decomposition of $V$ shows that $V$ is unitarily equivalent to $W$ or $W \oplus S^n$ or $W \oplus S^{*n}$ where $W$ is a unitary operator and $S$ is a unilateral shift of multiplicity one and then
\[
emac(F) = emac(V) = mac(W) + \frac {1}{2} |\mbox{ind } V|
\]
since $n = |\mbox{ind } V| = |\mbox{ind } F|$. Indeed $mac(W) = mac(\sJ W^*\sJ)$ while $emac(S^n \oplus S^{*n}) = n$ since $S^n \oplus S^{*n}$ is a finite rank perturbation of a bilateral shift of multiplicity $n$ for which $mac$ is $n$. We record the results of this discussion as the next Lemma.

\bigskip
\noindent
{\bf Lemma 4.1.} {\em
The essential multiplicity function $emac(F)$ of a {\rm mod} trace-class unitary operator is defined Haar a.e.\ on ${\mathbb T}$ and takes values in $\{0,1/2,1,3/2,2,\dots\infty\}$ and has the following properties
}

a) {\em If $F$ is unitary then $emac(F) = mac(F)$.
}

b) {\em If $F_1,f_2$ are {\rm mod} trace-class unitary and $W$ is a unitary operator so that $WF_1 - F_2W \in \sC_1$ then $emac(F_1) = emac(F_2)$.
}

c) {\em $emac(F_1 \oplus F_2) = emac(F_1) + emac(F_2)$.
}

d) {\em We have 
\[
emac(F^*)(e^{i\theta}) = emac(F)(e^{-i\theta})
\]
and if $\sJ$ is a conjugate-linear antiunitary operator the also
\[
emac(\sJ F\sJ)(e^{i\theta}) = emac(F)(e^{-i\theta}).
\]
}

e) {\em If $\mbox{ind } F \equiv 0(\mbox{\rm mod } 2)$ then $emac(F)$ takes values in $\{0,1,2,\dots,\infty\}$ and if $\mbox{ind } F \equiv 1(\mbox{mod } 2)$ then $emac(F)$ takes values in $\{1/2,3/2,5/2,\dots,\infty\}$. Moreover $emac(F) \ge 1/2|\mbox{ind } F|$.
}

f) {\em If $S$ is the unilateral shift operator in $\ell^2({\mathbb N})$ then
\[
emac(S^n) = emac(S^{*n}) = n/2.
\]
}

g) {\em $emac(F) = emac(pol(F))$}.

\bigskip
We will also need to prove a second lemma.

\bigskip
\noindent
{\bf Lemma 4.2.} {\em Let $F_1,F_2 \in \sU\sC_1(\sH)$ be so that $\sigma(F_j) \supset {\mathbb T}$, $j = 1,2$, $\mbox{ind }F_1 = \mbox{ind } F_2$ and $emac(F_1) = emac(F_2)$. Then there is a unitary operator $W$ so that $WF_1 - F_2W \in \sC_1$.
}

\bigskip
\noindent
{\em Proof.} To begin observe that the statement for the $F_j$'s is equivalent to that for the $F^*_j$'s, so we may assume $\mbox{ind } F_j \le 0$. Observe also that the spectrum condition is equivalent to $\sigma_e(F_j) = {\mathbb T}$. To prove the lemma we replace successively the $F_j$'s by others which are unitarily equivalent $\mbox{mod } \sC_1$ to them, and hence satisfy the same assumptions, till we get to an obvious assertion.  First we pass from $F_j$ to $pol(F_j)$ and then after a finite rank perturbation, the index being $\le 0$, we can assume $F_j$ is an isometry. Using then, for instance, the adapted Voiculescu theorem \cite{10} we can replace $F_j$ by $F_j \oplus G_j$ where $G_j$ is a unitary operator with singular spectrum. Passing to the Wold decomposition of the isometry $F_j$ we arrive at $E_j \oplus S^n \oplus G_j$ where $S$ is the unilateral shift and $E_1,E_2$ are unitary operators with equal $mac(E_j) = emac(F_j) - n/2$, $j = 1,2$. Choosing $G_1 = E_2 \mid \sH_{sing}(E_2)$, $G_2 = E_1 \mid \sH_{sing}(E_1)$ we arrive at unitary equivalence.\hfill\qed

\bigskip
Like in the case of a Hermitian operator $T$ in \cite{11}, where we defined an ordered group $\sF(T)$ and constructed an isomorphism with $K_0(\sE(T;\sC_1))$, we shall define an ordered group $\sF(U)$ and construct an isomorphism with $K_0(\sE(U;\sC_1))$, only $\sF(U)$ will be slightly more complicated to describe than $\sF(T)$. The elements of $\sF(U)$ are pairs $(f,u)$ where $f: {\mathbb T}\backslash\omega(U) \to 1/2 {\mathbb Z}$ is a measurable function up to almost everywhere equivalence, $n \in {\mathbb Z}$ and so that $|n| + |f(z)| \le C\ mac(U)(z)$ for some constant $C$ and $2f(z) \equiv n(\mbox{mod } 2)$ for almost all $z \in {\mathbb T}\backslash\omega(U)$. The operation on $\sF(U)$ is componentwise addition. Moreover $\sF(U)$ is an ordered group, the semigroup of elements which are $\ge 0$ being $\sF_+(U) = \{(f,n) \in \sF(U) \mid f(z) \ge |n|/2$ for $z \in {\mathbb T}\backslash\omega(U)\mbox{a.e.}\}$ Remark that unless $mac(U) \ge 1$, we must have $n = 0$ and $\sF(U)$ looks like the $\sF(T)$.

\bigskip
\noindent
{\bf Theorem 4.2.} {\em Let $U$ be a unitary operator with $\sigma(U) = \sigma_e(U)$. If $P \in \sP_h(\sE(U \otimes I_n;\sC_1))$ then $P(U \otimes I_n) \mid P\sH^n \in \sU\sC_1(P\sH^n)$. There exists a unique isomorphism
\[
femac(U): K_0(\sE(U;\sC_1)) \to \sF(U)
\]
so that if $P \in \sP_n(\sE(U \otimes I_n;\sC_1))$ then
\[
femac(U)([P]_0) = (emac(P(U \otimes I_n) \mid P\sH^n) \mid {\mathbb T}\backslash\omega(U),\mbox{ind } P(U \otimes I_n) \mid P\sH^n).
\]
}

\bigskip
\noindent
{\em Proof.} If $\sigma(U) \ne {\mathbb T}$ the remarks preceding the theorem show that the statement is equivalent to the result in \cite{11}. Thus we may assume $\sigma(U) = {\mathbb T}$. Since $[P,U \otimes I_n] \in \sC_1$ it follows that $P(U \otimes I_n) \mid P\sH^n \in \sU\sC_1(P\sH^n)$. To check that $femac(U)([P]_0)$ is well-defined let us first see that the formula depends only on $[P]_0$. Indeed if $[P]_0 = [Q]_0$ then $P \oplus I \otimes I_m \sim Q \oplus I \otimes I_m$ for some $m$. Then $(P \oplus I \otimes I_m)(U \otimes I_{n+m}) \mid (P\sH^n \oplus \sH^m)$ and $(Q \oplus I \otimes I_m)(U \otimes I_{n+m}) \mid (Q\sH^n \oplus \sH^m)$ are unitarily equivalent 
$\mbox{mod } \sC_1$ and this gives
\[
\begin{aligned}
&emac((P \oplus I \otimes I_m)(U \otimes I_{n+m}) \mid (P\sH^n \oplus \sH^m)) \\
&= emac((Q \oplus I \otimes I_m)(U \otimes I_{n+m}) \mid (Q\sH^n \oplus \sH^m))
\end{aligned}
\]
which gives
\[
emac(PU \mid P\sH^n) + memac(U) = emac(QU \mid Q\sH^n) + memac(U)
\]
and hence
\[
\begin{aligned}
&emac(PU \mid P\sH^n) \mid ({\mathbb T}\backslash\omega(U)) \\
&= emac(QU \mid Q\sH^n) \mid ({\mathbb T}\backslash\omega(U)).
\end{aligned}
\]
We also have clearly
\[
\begin{aligned}
&ind((P \oplus I \otimes I_m)(U \otimes I_{n+m}) \mid (P\sH^n \oplus \sH^m)) \\
&= ind((Q \oplus I \otimes I_m)(U \otimes I_{n+m}) \mid (Q\sH^n \oplus \sH^m)
\end{aligned}
\]
which means that
\[
ind(P(U \otimes I_n) \mid P\sH^n) = ind(Q(U \otimes I_n) \mid Q\sH^n).
\]
Next we need also to check that
\[
emac(P(U \otimes I_n) \mid P\sH^n) \mid ({\mathbb T}-\omega(U)),ind(P(U \otimes I_m) \mid P\sH^n)) \in \sF(U).
\]
We have $emac(P(U \otimes I_n) \mid P\sH^n) \ge 1/2 |ind(P(U \otimes I_n) \mid P\sH^n)|$ and $2emac(P(U \otimes I_n) \mid P\sH^n) \equiv ind(P(U \otimes I_n) \mid P\sH^n)(\mbox{mod } 2)$ by Lemma~4.1~e) applied to $F = P(U \otimes I_n) \mid P\sH^n$. If $G = (I \otimes I_n - P)(U \otimes I_n) \mid (I \otimes I_n - P)\sH^n$), then $F \oplus G$ is unitarily equivalent $\mbox{mod } \sC_1$ to $U \otimes I_n$ so that by Lemma~4.1~c), b) and a) we have
\[
emac(F) \le emac(F \oplus G) = emac(U \otimes I_n) = nmac(U).
\]
Thus $femac(U)([P]_0)$ is well-defined. Also since for $P \in \sP_h(\sE(U \otimes I_n;\sC_1))$, $Q \in \sP_h(\sE(U \otimes I_n;\sC_1)$ we have
\[
(P \oplus Q)(U \otimes I_{n+m}) \mid (P \oplus Q)\sH^{n+m} - (P(U \otimes I_n)P\sH^n) \oplus (Q(U \otimes I_m) \mid Q\sH^m) \in \sC_1
\]
using Lemma~4.1 we get
\[
femac(U)([P]_0) + femac(U)([Q]_0) = femac(U)([P \oplus Q]_0).
\]
This implies that $femac(U)$ extends to a homomorphism, which is also unique, of $K_0(\sE(U;\sC_1))$ into $\sF(U)$. We shall denote in the rest of the proof this extension still by $femac(U)$.

The next step is to check that $femac(U)$ is one-to-one and onto.

To check that $femac(U)$ is onto, it suffices to show that for every $(f,n) \in \sF_+(U)$ there is $P \in 
\sP_h(\sE(U \otimes I_m;\sC_1))$ for some $m \in {\mathbb N}$, so that $femac(U)([P]_0) = (f,n)$. Indeed, $(g,n) \in \sF(U)$ implies $(|g|,|n|) \in \sF(U)$ so $(g,n) = (|g|,|n|) - (|g|-g,|n|-n) \in \sF_+(U) - \sF_+(U)$. If $(f,n) \in \sF_+(U)$ let ${\tilde f}: {\mathbb T} \to \{0,1,2,\dots\}$ be defined as ${\tilde f}(z) = f(z) - |n|/2$ for $z \in {\mathbb T}\backslash\omega(U)$ and ${\tilde f}(z) = 0$ for $z \in \omega(U)$. Since ${\tilde f} \le C\ mac(U)$ we can find a projection $Q \in (U \otimes I_m)'$ for some $m$ so that $mac Q(U \otimes I_m) \mid Q\sH^m = {\tilde f}$. Thus $({\tilde f},0) = femac(U)([Q]_0)$ and we must find a projection $P$ so that $(|n|/2,n) = femac(U)([T]_0)$. If $n \ne 0$, $(|n|/2,n) \in \sF(U)$ implies that $mac(U) \ge 1$ on ${\mathbb T}$. Then there is a projection ${\tilde P} \in \sP_h(\sE(U;\sC_1))$, ${\tilde P} \in (U)'$ so that $mac({\tilde P}U \mid {\tilde P}\sH) = 1$. Thus ${\tilde P}U \mid {\tilde P}\sH$ is unitarily equivalent to the bilateral shift in $\ell^2({\mathbb Z})$ and hence we can find a projection ${\overset{\approx}{P}} \in \sP_h(\sE(U;\sC_1))$ so that ${\overset {\approx}{P}}U \mid {\overset{\approx}{P}}\sH$ is unitarily equivalent to the unilateral shift of multiplicity $1$ if $n > 0$. It follows that if $P = {\overset{\approx}{P}} \otimes I_{|n|}$ we will have
\[
femac(U)([P]_0) = (|n|/2,n).
\]
Note that this also proves that
\[
\sF_+(U) \subset femac(U)((K_0(\sE(U;\sC_1))_+).
\]

To show $femac(U)$ is one-to-one it suffices to show that if $P,Q \in \sP_h(\sE(U \otimes I_n,\sC_1))$ are so that $femac(U)([P]_0) = femac(U)([Q]_0)$ then $I \oplus P \sim I \oplus Q$. The equality of $femac(U)$ for $[P]_0$ and $[Q]_0$ gives that
\[
ind(P(U \otimes I_n) \mid P\sH^n) = ind(Q(U \otimes I_n) \mid Q\sH^n)
\]
and
\[
emac(P(U \otimes I_n) \mid P\sH^n) \mid {\mathbb T}\backslash\omega(U) = emac(Q(U \otimes I_n) \mid Q\sH^n) \mid {\mathbb T}\backslash\omega(U).
\]
Note that we can replace $P,Q$ by $P \oplus I$, $Q \oplus I$ and $n$ by $n+1$ and the Fredholm indices don't change, but the equality of the $emac$ will extend to all of ${\mathbb T}$. Then 
\[
F_1 = (P \oplus I)(U \otimes I_{n+1}) \mid (P \oplus I)\sH^{n+1}
\]
and
\[
F_2 = (Q \oplus I)(U \otimes I_{n+1}) \mid (Q \oplus I)\sH^{n+1}
\]
satisfy the assumptions of Lemma~$4.2$. Then the unitary equivalence $\mbox{mod } \sC_1$ of $F_1$ and $F_2$ which the lemma asserts, gives $P \oplus I \sim Q \oplus I$.

Thus $femac(U)$ is a bijection and to conclude the proof we need only to show that $femac(U)([P]_0) \in \sF_+(U)$. This follows from Lemma~$4.1$~e) applied to $P(U \otimes I_n) \mid P\sH^n$.\hfill\qed

\section{A $K$-theory problem}
\label{sec5}

We point out in this section the simplest case of the problem of computing the $K_1$-group, accompanied by a few remarks.

\bigskip
\noindent
{\bf Problem.} Let $T$ be the Hermitian operator of multiplication by the coordinate function in $L^2([0,1],d\lambda)$, $(Tf)(x) = xf(x)$, $d\lambda$ Lebesgue measure. What is $K_1(\sE(T,\sC_1)$? Is it non-trivial?

\bigskip
An example of a unitary operator $U$ in $\sE(T,\sC_1)$ about the triviality of the $K_1$-class of which one may wonder is the following. Using for instance the absorption version of the Voiculescu theorem \cite{10}, there is a unitary operator
\[
V: L^2([0,1],d\lambda) \to L^2([0,1],d\lambda) \oplus \ell^2({\mathbb N}) \oplus \ell^2({\mathbb N})
\]
so that
\[
VTV^* - T \oplus \alpha I \oplus \beta I \in \sC_1
\]
where $0 \le \alpha < \beta \le 1$. Clearly this shows that $\sE(T;\sC_1)$ and $\sE(T \oplus \alpha I \oplus \beta I)$ are isomorphic. Let $S$ be the unilateral shift operator in $\ell^2({\mathbb N})$ and let
\[
U = R + I \oplus S \oplus S^*
\]
where $R$ is a rank one partial isometry from $0 \oplus 0 \oplus \ker S^*$ to $0 \oplus \ker S^* \oplus 0$ so that $I \oplus S \oplus S^* + R$ is unitarily equivalent to $I \oplus W$ where $W$ is a bilateral shift operator i.e.\ a unitary operator. Is $[U]_1$ trivial or not?

A general remark one may find useful is that by \cite{11} we have $K_0(\sK(T,\sC_1)) = K_0(\sK) = {\mathbb Z}$ and $K_1(\sK(T,\sC_1)) = K_1(\sK) = 0$ so that the 6-terms $K$-theory exact sequence for $0 \to \sK(T,\sC_1) \to \sE(T,\sC_1) \to \sE/\sK(T,\sC_1) \to 0$ gives that $K_1(\sE(T,\sC_1))$ is isomorphic to $\ker \partial \subset K_1(\sE/\sK(T,\sC_1))$ where $\partial$ is the connecting map given by the Fredholm index from $K_1(\sE/\sK(T;\sC_1)) \to K_0(\sE(T;\sC_1))$ (which is easily seen to be surjective). Recall also that $\sE/\sK(T;\sC_1)$ is isomorphic to a $C^*$-algebra with center generated by the class of $T$ in $\sB/\sK$ i.e.\ $T + \sK$. The joint spectrum of the classes of $T$ and $U$ in $\sE/\sK(T;\sC_1)$ is then
\[
([0,1] \times \{1\}) \cup (\{\alpha,\beta\} \times {\mathbb T}) \subset {\mathbb C} \times {\mathbb C}
\]
where ${\mathbb T} = \{z \in {\mathbb C} \mid |z| = 1\}$. This also shows that $\sE/\sK(T;\sC_1)$ is far from being the $C^*$-algebra of a continuous field of $C^*$-algebras over $[0,1]$ the spectrum of its center.

We recall that if $B$ and $A$ are Hermitian operators without singular spectrum and $B - A \in \sC_1$, then the Kato--Rosenblum theorem (\cite{6}, \cite{7}, \cite{8}) implies that the strong limit
\[
W^+(B,A) = s - \lim_{t \to +\infty} e^{itB}e^{-itA}
\]
exists and is a unitary operator intertwining $B$ and $A$
\[
BW^+(B,A) = W^+(B,A)A.
\]
A consequence of this is the following result about $K_1(\sE(T;\sC_1))$.

\bigskip
\noindent
{\bf Proposition 5.1}. {\em Let $U \in \sE(T \otimes I_n;\sC_1) \simeq \sM_n(\sE(T;\sC_1))$ be a unitary operator. Then in $K_1(\sE(T;\sC_1))$ we have
\[
[U]_1 = [W^+(U(T \otimes I_n)U^*,T \otimes I_n)]_1.
\]
}

\bigskip
\noindent
{\em Proof.} Since
\[
U(T \otimes I_n)U^* = W^+(U(T \otimes I_n)U^*,T \otimes I_n)(T \otimes I_n)(W^+(U(T \otimes I_n)U^*,T \otimes I_n))^*
\]
we have
\[
U^*W^+(U(T \otimes I_n)U^*,T \otimes I_n) \in (T \otimes I_n)'.
\]
The commutant $(T \otimes I_n)'$ is a von~Neumann algebra so $K_1((T \otimes I_n)') = 0$ and this implies
\[
[U^*W^+(U(T \otimes I_n)U^*,T \otimes I_n)]_1 = 0
\]
which is the desired result.\hfill\qed

\bigskip
The preceding proposition used only the fact that the spectrum of $T$ is absolutely continuous. Also, the proposition holds with $W^+$ replaced by $W^-(B,A) = W^+(-B,-A)$.

The interest of Proposition~$5.1$ for the $K_1$-problem is that the $K_1$-classes are precisely the classes of $W^+(A,T \otimes I_n)$ where $A-T \otimes I_n \in \sC_1$ is so that $W^+(A,T \otimes I_n)$ is unitary. Thus the problem can be rephrased in terms of special trace-class perturbations of the $T \otimes I_n$.


\begin{thebibliography}{99}

\bibitem{1} Birman, M.~S., and Krein, M.~G., {\em On the theory of wave operators and scattering operators}, Dokl.  Akad. Nauk SSSR {\bf 144}, 475--478 (1962) (Sov. Math. Dokl. 3, 740--747 (1962).

\bibitem{2} Blackadar, B., {\em $K$-theory for operator algebras}, MSRI Publications, Vol.~5, Springer Verlag, 1986.

\bibitem{3} Bourgain, J., and Voiculescu, D.~V., {\em The essential centre of the {\rm mod}-a-diagonalization ideal commutant of an $n$-tuple of commuting Hermitian operators}, Non-commutative analysis, operator theory and applications, 77--80, Oper. Theory Adv. Appl. {\bf 252}, Linear Oper. Linear Syst., Birkhauser/Springer (2016).

\bibitem{4} Gohberg, I.~C., and Krein, M.~G, {\em Introduction to the theory of linear non-selfadjoint operators}, Translations of Mathematical Monographs, Vol.~18, AMS, Providence, RI (1969).

\bibitem{5} Johnson, B.~E., and Parrott, S.~K., {\em Operators commuting modulo the set of compact operators with a von~Neumann algebra}, J. Funct. Anal. {\bf 11} (1972), 39--61.

\bibitem{6} Kato, T., ``Perturbation Theory for Linear Operators'', Classics of Mathematics, Springer Verlag, 1995.

\bibitem{7} Putnam, C.~R., ``Commutation Properties of Hilbert Space Operators and Related Topics'', Springer Verlag, Berlin--Heidelberg--New~York, 1967.

\bibitem{8} Reed, M., and Simon, B., {\em Methods of modern physics}, Vol.~III: Scattering theory, Academic Press, 1979.

\bibitem{9} Simon, B., {\em Trace ideals and their applications}, 2nd Ed., Mathematical Surveys and Monographs, Vol.~120, AMS, Providence, RI, 2005.

\bibitem{10} Voiculescu, D.~V., {\em Some results on norm-ideal perturbations of Hilbert space operators I}, J.~Operator Theory {\bf 1} (1979), 3--37.

\bibitem{11} Voiculescu, D.~V., {\em $K$-theory and perturbations of absolutely continuous spectra}, Comm. Math. Phys. {\bf 365} (2019), No.~1, 363--373.

\bibitem{12} Voiculescu, D.~V., {\em Commutants {\em mod} normed ideals}, arXiv:1810.12497.

\end{thebibliography}
\end{document}